\magnification=\magstep0
\input amstex
\define\la{\lambda}
\documentstyle{amsppt}
\NoBlackBoxes
\NoRunningHeads
\nologo
\topmatter
\title{Combinatorial Orthogonal Expansions}
\endtitle
\author A. de M\'edicis 
\footnote{\hbox{This work was supported by NSERC funds.}\hfill}
and D. Stanton
\footnote{\hbox{This work was supported by NSF grant DMS-9001195.\hfill}}
\endauthor
\address School of Mathematics, University of Minnesota,
Minneapolis, MN 55455.
\endaddress
\address School of Mathematics, University of Minnesota,
Minneapolis, MN 55455.
\endaddress
\abstract{The linearization coefficients for a set of 
orthogonal polynomials are given explicitly as a weighted sum of combinatorial
objects. Positivity theorems of Askey and Szwarc are 
corollaries of these expansions.}  
\endabstract
\endtopmatter

\subheading{1. Introduction}
Given a set of orthogonal polynomials $p_n(x)$, the linearization coefficients
$a_{mn}^k$ are given by 
$$ 
p_m(x)p_n(x)=\sum_k a_{mn}^k p_k(x).
$$
Askey \cite{1} and Szwarc \cite{4,5} have given sufficient conditions 
on the three-term recurrence relation coefficients $\alpha_n$, $\beta_n$, 
and $\gamma_n$ in
$$
\alpha_{n+1}p_{n+1}(x)=(x-\beta_n)p_n(x)-\gamma_{n-1} p_{n-1}(x)
\tag1.1
$$ 
so that
$a_{mn}^k$ is non-negative. In this paper we give in Theorem 1 and Theorem 2
explicit formulas for $a_{nm}^k$ as a polynomial in the 
$\alpha_j's$, $\beta_j's$ and the $\gamma_j's$, which give these theorems.

The idea is to represent $a_{mn}^k$ as a generating function of paths, whose
weights are products of differences. Monotonicity hypotheses on the 
coefficients force the weights to be individually positive, these are the conditions in \cite{1} and \cite{4}. For example, if $p_n(x)$ is monic; 
$\alpha_n=1$, $\beta_n=b_n$, and $\gamma_{n}=\la_{n+1}$, we have
$$
\align
a_{33}^3=&
(b_3-b_0)(b_3-b_1)(b_3-b_2)+(b_3-b_0)\la_4+(b_3-b_0)(\la_3-\la_2)+\\
&(b_4-b_1)\la_4+(b_3-b_2)\la_4+(b_2-b_1)\la_3+(b_3-b_2)(\la_3-\la_1).
\tag1.2
\endalign
$$
If $b_j$ and $\la_j>0$ are increasing, then $a_{33}^3$ is non-negative, see
\cite{1}.

\subheading{2. The theorems} We first recall some terminology 
and results in \cite{3} and \cite{6}.

We let $L$ denote the positive definite linear functional on the space of polynomials which corresponds to the orthogonal polynomials (1.1). So 
$L(x^n)=\mu_n$, the $nth$ moment of a measure for $p_n(x)$. It is easy to see
that
$$
a_{mn}^k=L(p_mp_np_k)/L(p_kp_k).
$$
Since $L(p_kp_k)=\gamma_{0}\cdots\gamma_{k-1}/\alpha_1\cdots\alpha_k>0$, 
we find instead $L(p_mp_np_k)$.

Viennot \cite{6} gave a combinatorial interpretation for the polynomials
$p_n(x)$ and their moments $\mu_n$, in terms of pavings and Motzkin paths respectively. We review these terms below.

A {\it{Motzkin path}} $P$ is a lattice path in the plane, which lies at or 
above the $x$-axis, and has steps of $(1,0)$ (horizontal=$H$), 
$(1,1)$ (up=$U$), or $(1,-1)$ (down=$D$).
The weight of a path $P$, $w(P)$, is defined by the product of the 
weights of its individual edges,
$$
w(P)=\prod_{edges\thinspace e} w(e).
\tag2.1
$$

A {\it{paving}} $\pi$ of the integers $\{1,\cdots,k\}$ is a 
collection of disjoint sets of cardinalities 1 (called monominos), 
and 2 (called dominos). 
The elements of a domino must be consecutive integers. 
For example, 
$\{ \{2,3\}, \{5\}, \{6,7\},\{9\}\}$ is a paving of $\{1,\cdots,9\}$. 
Points not in any of the sets are called {\it{isolated}}. The weight of 
a paving is defined to be the product of the individual 
weights of the monominos, dominos, and isolated points.

For Askey's theorem we need a special weight on edges $e$ of a Motzkin path. Suppose the path $P$ begins at $(0,m)$ 
and ends at $(k,n)$. We define
$$
w(\text{edge starting at }(i,j))=
\cases 
(b_j-b_i) \text{ if the edge is H,}\\
(\la_j-\la_{i+1})\text{ if the edge is D, and followed by U,}\\
\la_j\text{ if the edge is D, and not followed by U,}\\
1\text{ if the edge is U.}\\
\endcases
\tag2.2
$$

\proclaim{Theorem 1} Suppose that $\alpha_n=1$, $\beta_n=b_n$, and $\gamma_{n}=\la_{n+1}$. Then
$$
L(p_mp_np_k)=\la_{1}\cdots\la_n\sum_{P} w(P),
$$
where $P$ is a Motzkin path from $(0,m)$ to $(k,n)$,
and $w(P)$ is given by (2.1) and (2.2).
\endproclaim

For example, if $k=m=n=3$ in Theorem 1, there are 7 Motzkin paths
from (0,3) to (3,3): $HHH$, $HUD$, $HDU$, $UHD$, $UDH$, $DHU$, $DUH$. 
The weights of these 7 paths are the 7 terms in (1.2).

\demo{Proof of Theorem 1} One can prove that both sides in Theorem 1
have the same recurrence relation, which is given in \cite{1}. 

An alternative proof is to use Viennot's combinatorial 
interpretation for \newline $L(p_mp_np_k)/\la_1\cdots\la_n$, \cite{6}. 
It is the generating function for ordered pairs 
$(P,\pi)$, where $P$ is a Motzkin path from 
$(0,m)$ to $(l,n)$, and $\pi$ is a paving of the integers $\{1,\cdots,k\}$  
with $l$ isolated integers. The weight of $(P,\pi)$ is the product 
of the weights of $P$ and $\pi$. In $P$, an up edge starting at $(i,j)$ 
has weight $1$, a down edge $\lambda_j$, and an across edge $b_j$.
For $\pi$, a monomino at $\{i\}$ has weight $-b_{i-1}$, 
and a domino at $\{i,i+1\}$ has weight $-\la_i$.  

Given $(P,\pi)$ we create a unique path $P'$ by inserting in $P$, as
the $ith$ step of $P'$, an $H$ edge 
if $\pi$ has a monomino in position $i$. If 
$\pi$ has a domino starting in position $i$, we insert two steps, $DU$, 
in $P$, for the $ith$ and $(i+1)st$ steps of $P'$. 
The result is a single path $P'$ from 
$(0,m)$ to $(k,n)$. The weight of the path is given by (2.2): 
the negative terms correspond to the weight in $\pi$, 
the positive terms to the weight in $P$. 
\qed\enddemo

It is easy to see that Theorem 1 implies Askey's theorem.

\proclaim{Corollary 1} If $\la_j$ and $b_j$ are increasing, 
with $\la_j>0$, then $a_{mn}^k\ge 0$.
\endproclaim
\demo{Proof} 
We can assume by symmetry that $k\le n$, 
Then it is clear that each vertex $(i,j)$ in $P$ satisfies $i\le j$. 
Thus all weights are non-negative if the $b_j$'s and $\la_j$'s are 
increasing.
\qed\enddemo

Theorem 1 can be restated in terms of walks of length $m$ on the 
non-negative integers, starting at $k$, and ending at $n$, 
with steps of size $+1$, $-1$, or $0$. 

We let $p_n'(x)$ be another set of orthogonal polynomials satisfying 
$$
\alpha_{n+1}'p_{n+1}'(x)=(x-\beta_n')p_n'(x)-\gamma_{n-1}'p_{n-1}'(x).
$$ 
More generally, we consider
$$
p_m(x)p_k'(x)=\sum_{n} b_{mk}^n p_n(x).
\tag2.3
$$

It is clear that $b_{mk}^n=L(p_mp_k'p_n)/L(p_np_n)$. 
We will give an interpretation for $L(p_mp_k'p_n)$, which is 
non-negative when $b_{mk}^n$ is, since $L$ is positive definite.

We generalize Szwarc's theorem by finding a combinatorial 
interpretation for $L(p_mp_k'p_n)$ in (2.3). A {\it{generalized Motzkin path}}
allows a fourth type of edge: HH (across by two units). We define a weight 
$v(P)$ on 
generalized Motzkin paths from $(0,m)$ to $(k,n)$ again as a product of 
weights of edges, 
$$
v(\text{edge starting at }(i,j))=
\cases 
(\beta_j-\beta_{i}') \text{ if the edge is H,}\\
(\gamma_j-\alpha_{i}')\text{ if the edge is U, and preceded by D,}\\
\gamma_j\text{ if the edge is U, and not preceded by D,}\\
(\alpha_j-\alpha_i')\text{ if the edge is D, and preceded by U,}\\
\alpha_j\text{ if the edge is D, and not preceded by U,}\\
(\alpha_j+\gamma_j-\alpha_i'-\gamma_i')\alpha_{i+1}'\text{ if the edge is HH, preceded by U or D,}\\
(\alpha_j+\gamma_j-\gamma_i')\alpha_{i+1}'\text{ if the edge is HH, not 
preceded by U or D.}\\
\endcases
\tag2.4
$$

\proclaim{Theorem 2} We have 
$$
L(p_mp_np_k')=\frac{\gamma_0\cdots \gamma_{k-1}}
{\alpha_1\cdots\alpha_m\alpha_1'\cdots\alpha_k'}\sum_{P} v(P),
$$
where $P$ is a generalized Motzkin path from $(0,m)$ to $(k,n)$,
and $v(P)$ is given by (2.1) and (2.4).
\endproclaim

\demo{Proof} Again we will use Viennot's interpretation for 
$L(p_mp_np_k')\alpha_1\cdots\alpha_m/\gamma_0\cdots \gamma_{k-1}$. 
The weights on the edges, monominos, and dominos slightly change. 
Let $P'$ denote the Motzkin path and $\pi'$ the paving. 
In $P'$, the $U$, $D$, $H$ edges starting at $(i,j)$ have weights
 $\gamma_j$, $\alpha_j$, and $\beta_j$ respectively. In $\pi'$, a monomino
$\{i\}$ has weight $-\beta_{i-1}'/\alpha_i'$, a domino $\{i,i+1\}$ 
has weight $-\gamma_{i-1}'\alpha_i'/(\alpha_i'\alpha_{i+1}')$, 
and an isolated point $i$ has weight $1/\alpha_i'$. 
Note that every paving has a factor 
of $1/\alpha_1'\cdots\alpha_k'$. We therefore disregard the 
denominators of the weights of the pavings, and put this constant 
factor in the statement of Theorem 2. 

As in Theorem 1, we will merge pavings $\pi'$ with the paths $P'$ 
to create a generalized Motzkin path $P$ whose weights are given by 
(2.1) and (2.5) 
$$
u(\text{edge starting at }(i,j))=
\cases 
(\beta_j-\beta_{i}') \text{ if the edge is H,}\\
\gamma_j\text{ if the edge is U,}\\
\alpha_j\text{ if the edge is D,}\\
-\gamma_i'\alpha_{i+1}'\text{ if the edge is HH.}\\
\endcases
\tag2.5
$$
The basic idea is to insert certain edges into $P'$ to create $P$, while 
simultaneously deleting all monominos and dominos in $\pi'$. This is done by 
inserting an $H$ edge in $P'$ starting at $(i,j)$, if $\pi'$ has the monomino 
 $\{i+1\}$. We insert an $HH$ edge in $P'$ starting at $(i,j)$, if $\pi'$ has the domino $\{i+1,i+2\}$. We obtain a multiset of generalized Motzkin paths 
$P:(0,m)\rightarrow (k,n)$, from which the multiplicities are eliminated 
by using the weight (2.5).

Let $S$ be the set of all generalized Motzkin paths from $(0,m)$ to $(k,n)$.
We just found that the linearization coefficients are, up to a constant, the generating function for $S$ with weight (2.5). We want weight (2.4) 
instead of (2.5). We will do this via an involution.

The (2.4) weights of the edges of $P\in S$ are not monomials, 
instead they are sums of monomials.   
Thus we can consider the multiset $M_1$ of 
paths $P\in S$, where 
the multiplicity of $P$ in $M_1$ 
is the product of the number of monomials in the weight of 
the edges $e\ne H$ of $P$. The weight of any element of $M_1$ is the 
product of a choice of monomials for each edge.
On $M_1$ we will construct a weight-preserving 
sign-reversing involution, whose fixed point set consists of 
all paths $P$ exactly once, with weights (2.5).

It remains to give the involution $\Phi$ on the multiset $M_1$ of paths $P$. 
Note that we want 
to eliminate all weights in the edges that include $\alpha'$, 
except for the $-\gamma_i'\alpha_{i+1}'$ term in $HH$. Scan the 
path $P$ from right to left, and find the first such term in the choice of 
monomials for the weights. Suppose the edge containing this term is $HH$, 
preceded by 
$U$ or $D$. From (2.5), the weight we need to eliminate is one term from 
$(\alpha_j+\gamma_j-\alpha_i')\alpha_{i+1}'$. If the preceding edge is $D$, 
replacing the $HH$ edge by a pair $UD$ will cancel the  $(\gamma_j-\alpha_i')\alpha_{i+1}'$ terms, while replacing the $HH$ 
edge by $DU$ will cancel the $\alpha_j\alpha_{i+1}'$ term. Similarly, 
if the preceding edge to $HH$ is $U$, replacing $HH$ by $UD$ and $DU$ will
cancel the $\gamma_j\alpha_{i+1}'$ and 
$(\alpha_j-\alpha_i')\alpha_{i+1}'$ terms, respectively. 
If the first edge containing $\alpha'$ is $HH$, not preceded by $U$ or $D$,
we must eliminate  $(\alpha_j+\gamma_j)\alpha_{i+1}'$. This time 
replacing $HH$ by $DU$ and $UD$ eliminates a single term each. 

This defines $\Phi(P)=Q$, when the first appropriate $\alpha'$ edge of $P$ is $HH$.
If the first appropriate $\alpha'$ edge of $P$ is not $HH$, then $\alpha'$
 must be a choice of weight from a $DU$ or $UD$. Then we invert the 
previous case.
It is easy to check that the involution $\Phi$ is 
well defined on $M_1$, with the 
stated fixed points.
\qed\enddemo

Corollary 2 generalizes \cite{4, Theorem 2}.

\proclaim{Corollary 2} If $\alpha_i,\alpha_i',\gamma_i,\gamma_i'>0$, 
$\beta_j\ge \beta_i'$,  
$\alpha_j\ge \alpha_i'$, 
$\alpha_j+\gamma_j\ge \alpha_i'+\gamma_i'$,
$\gamma_j\ge \alpha_i'$,
for $j\ge i$, and $k\le {\text{max}}\{m,n\}$, then $b_{mk}^n\ge 0$.
\endproclaim
\demo{Proof} 
Assume $k\le n$. 
The inequalities insure that the individual weights in Theorem 2 are 
positive, since the indices of the primed variables cannot be greater 
than the indices of the unprimed variables. By symmetry we obtain the 
$k\le{\text{max}}\{m,n\}$ case.
\qed\enddemo

The connection coefficient problem is the $m=0$ special case of Theorem 2.
Non-zero coefficients occur only for $k\ge n$. In this case, 
along our path $P$, vertices $(i,j)$ satisfy $i\ge j$, 
so we assume the inequalities of Corollary 2 
hold in this range.  This implies Askey's theorem in \cite{2}. 

The theorems in \cite{5} can also be generalized, for example:

\proclaim{Corollary 3} If $\beta_j=\beta_i'=0$,
$\alpha_i,\alpha_i',\gamma_i,\gamma_i'>0$,   
$\alpha_{2j}\ge \alpha_{2i}'$,
$\alpha_{2j+1}\ge \alpha_{2i+1}'$,
$\alpha_{2j}+\gamma_{2j}\ge \alpha_{2i}'+\gamma_{2i}'$,
$\alpha_{2j+1}+\gamma_{2j+1}\ge \alpha_{2i+1}'+\gamma_{2i+1}'$,
$\gamma_{2j}\ge \alpha_{2i}'$,
$\gamma_{2j+1}\ge \alpha_{2i+1}'$,
for $j\ge i$, $m$ is even, and $k\le n$, then $b_{mk}^n\ge 0$.
\endproclaim
\demo{Proof} Under the assumption that $m$ is even, and all $\beta's =0$,
all vertices $(i,j)$ on the path $P$ of Theorem 2 have the property 
that $i$ and $j$ have the same parity.
\qed\enddemo

\enddocument
\end